\documentclass[12pt]{article}
\usepackage{amsmath, amssymb, amsthm, graphicx, geometry}
\usepackage{hyperref} 

\title{Ill-Posed Configurations in Random and Experimental Data Points Collection}
\author{\Large{Netzer Moriya}}
\date{}

\begin{document}
\maketitle

\begin{abstract}
Ill-posed configurations, such as collinear or coplanar point arrangements, are a persistent challenge in computational geometry, 
complicating tasks as in triangulation and convex hull construction. This paper discusses the probability of such 
configurations arising in two scenarios: (1) data sampled randomly from a uniform distribution, and (2) data collected from physical 
systems, such as reflective surfaces or structured environments. 
We present a probabilistic framework, analyze the geometric and sampling constraints, and provide some mathematical insights into how data acquisition processes influence the likelihood of degeneracies.
Notably, our findings reveal that degeneracies occur more frequently in physical systems than in purely random simulations due to systematic biases introduced by instrumental setups and environmental structures, emphasizing the risks of drawing conclusions solely based on assumptions derived from random data.
\end{abstract}

\section{Introduction}

The stability and robustness of geometric computations often hinge on the nature of input data. Ill-posed configurations, 
such as collinear or coplanar arrangements, lead to numerical instabilities and pose significant challenges, particularly in 
problems like the Smallest Enclosing Sphere (SES)~\cite{badoiu2003coresets}, Convex Hull 
construction~\cite{convexHullSymmetry, convexHullParallel} , model facets reconstruction~\cite{delaunayReconstruction, 
featurePreservingReconstruction} and 
triangulation~\cite{learnableTriangulation}. In datasets sampled uniformly at random, 
the probability of encountering such degeneracies decreases with increasing dimensionality, as higher-dimensional spaces 
inherently reduce the likelihood of points aligning in lower-dimensional subspaces.

In contrast, physical data collection methods, such as LiDAR or photogrammetry, introduce systematic biases, quantization effects, digitization, measurement noise and environmental 
constraints, making degeneracies far more common. For instance, point clouds generated from 
Unoccupied Aerial Systems (UAS) equipped with LiDAR sensors can display inherent geometric regularities due to the structured 
nature of the scanned environments \cite{gomez2024uas}. These regularities can lead to collinear or coplanar point arrangements, 
complicating subsequent geometric computations.

Despite extensive research on degeneracies in both random and structured 
datasets~\cite{Bakushinsky1994, Burgisser2008}, a comprehensive framework systematically 
comparing the likelihood of ill-posed configurations across these contexts remains underdeveloped. This study seeks to 
address this gap by integrating probabilistic analyses with empirical observations derived from simulated scenarios that 
model physical systems, thereby providing deeper insights into how data acquisition processes influence the occurrence of degeneracies.

This paper analyzes degeneracies in two contexts: (1) statistical data sampled uniformly and independently from \( \mathbb{R}^d \), 
representing idealized random configurations, and (2) simulated data derived from scenarios that model physical systems,
influenced by simulated sensor biases and environmental structure. Purely statistical data provides a baseline for probabilistic analysis, while simulating data from models of physical systems reflects practical constraints, including noise and geometric regularities. 

In particular, degeneracies are more frequent in physical systems than in purely random simulations, primarily due to systematic biases introduced by sensor geometry and environmental structures. This highlights the risks of relying solely on conclusions derived from random data when analyzing real-world datasets.

This distinction is central to the comparative framework presented.

\section{Problem Statement}

We consider the probability of encountering ill-posed configurations in two contexts: 
\begin{enumerate}
    \item points sampled uniformly from \( [0,1]^d \subset \mathbb{R}^d \), where each point \( p_i \) is independently drawn from a uniform distribution \( U(0,1) \).
    \item data collected from simulated physical systems, where points \( p_i \) can be modeled as \( p_i = s_i + \epsilon_i \), with \( s_i \) lying on a structured.
surface \( S \subset \mathbb{R}^d \) and \( \epsilon_i \) representing noise.
\end{enumerate}
The first scenario analyzes intrinsic probabilistic properties, focusing on degeneracies arising purely from random configurations. 
In contrast, the second scenario incorporates systematic effects, including noise and biases from sensor geometry, as well as 
geometric regularities introduced by reflective surfaces in structured environments, such as lines, planes, and curved surfaces. 
These regularities greatly increase the likelihood of degeneracies compared to purely random data. This distinction enables a 
clearer linkage between theoretical randomness and real-world constraints in degeneracy modeling.

\subsection{Definitions}  
Let \( P = \{p_1, p_2, \ldots, p_N\} \subset \mathbb{R}^d \) denote a set of \( N \) points in \( d \)-dimensional Euclidean space. A subset \( P_k \subset P \) is said to be \textbf{degenerate} if there exists a proper subspace \( M \subset \mathbb{R}^d \) such that all points in \( P_k \) satisfy \( p_i \in M \) for \( p_i \in P_k \).

\begin{itemize}
    \item \textbf{Collinear Points:} A subset \( P_k \subset P \) is collinear if there exists a one-dimensional affine subspace \( M \subset \mathbb{R}^d \) (i.e., a line, \( \dim(M) = 1 \)) such that \( p_i \in M \) for all \( p_i \in P_k \).
    \item \textbf{Coplanar Points:} A subset \( P_k \subset P \) is coplanar if there exists a two-dimensional affine subspace \( M \subset \mathbb{R}^d \) (i.e., a plane, \( \dim(M) = 2 \)) such that \( p_i \in M \) for all \( p_i \in P_k \).
    \item \textbf{Nearly Spherical Points:} A subset \( P_k \subset P \) is nearly spherical if there exists a sphere \( S \subset \mathbb{R}^d \) with center \( c \in \mathbb{R}^d \) and radius \( r > 0 \) such that the points in \( P_k \) satisfy \( \|p_i - c\| \approx r \) for all \( p_i \in P_k \), within a given tolerance \( \epsilon > 0 \).
\end{itemize}

\subsection{Probabilistic Framework}  
The probability that a subset \( P_k \) is degenerate can be expressed as:
\[
P_k(\text{degenerate}) = \mathbb{P}(P_k \subseteq M),
\]
where \( M \) denotes a lower-dimensional subspace of \( \mathbb{R}^d \) satisfying the degeneracy condition.

For a set of \( N \) points, the probability of encountering at least one degenerate subset across all possible subsets is given by:
\[
P_{\text{degenerate}} = 1 - \prod_{k=1}^{N} \left(1 - P_k(\text{degenerate})\right)^{\binom{N}{k}},
\]
where \( \binom{N}{k} \) is the binomial coefficient, representing the number of subsets of size \( k \).

\section{Degeneracies in Purely Random Data}

Let \( P = \{p_1, p_2, \ldots, p_N\} \) represent \( N \) points sampled uniformly and independently from the unit hypercube \( [0,1]^d \subset \mathbb{R}^d \). A subset \( P_k \subset P \), where \( |P_k| = k \), is said to be degenerate if all points in \( P_k \) lie in a proper affine subspace of \( \mathbb{R}^d \). We analyze this degeneracy for random data, focusing on probabilities and scaling behavior, and provide explicit derivations for both collinear and coplanar configurations.

\subsection{Degeneracy Conditions\texorpdfstring{\protect\footnote{Hereafter we ignore degeneracies involving higher-dimensional affine subspaces, such as Nearly Spherical configurations due to their inherent complexities.}}{}}


\begin{itemize}
\item \textbf{Collinearity in \( \mathbb{R}^3 \):} For a subset of three points \( \{p_1, p_2, p_3\} \subset \mathbb{R}^3 \), collinearity implies that all three points lie on a single straight line. This is equivalent to the condition that the vectors formed by the points are linearly dependent. Let \( \mathbf{v}_1 = p_2 - p_1 \) and \( \mathbf{v}_2 = p_3 - p_1 \). The points are collinear if the cross product of these vectors is zero:
\[
\mathbf{v}_1 \times \mathbf{v}_2 = \mathbf{0}, \quad \text{where } \mathbf{v}_1 = 
\begin{bmatrix}
x_2 - x_1 \\
y_2 - y_1 \\
z_2 - z_1
\end{bmatrix}
\text{ and } \mathbf{v}_2 = 
\begin{bmatrix}
x_3 - x_1 \\
y_3 - y_1 \\
z_3 - z_1
\end{bmatrix}.
\]
Alternatively, this can be expressed in determinant form:
\[
\det\begin{bmatrix}
x_2 - x_1 & y_2 - y_1 & z_2 - z_1 \\
x_3 - x_1 & y_3 - y_1 & z_3 - z_1
\end{bmatrix} = 0.
\]
A zero cross product or determinant indicates that the points are collinear in \( \mathbb{R}^3 \).

    \item \textbf{Coplanarity in \( \mathbb{R}^3 \):} For a subset of four points \( \{p_1, p_2, p_3, p_4\} \subset \mathbb{R}^3 \), coplanarity implies that all four points lie on a single plane. This is equivalent to the determinant condition:
    \[
    \det\begin{bmatrix}
    x_1 & y_1 & z_1 & 1 \\
    x_2 & y_2 & z_2 & 1 \\
    x_3 & y_3 & z_3 & 1 \\
    x_4 & y_4 & z_4 & 1
    \end{bmatrix} = 0.
    \]
    This determinant tests whether the four points span a three-dimensional space; a zero determinant indicates that the points are coplanar.

    \item \textbf{Nearly Spherical points in \( \mathbb{R}^3 \):} While nearly spherical degeneracies are conceptually important, they are exceedingly rare in purely random datasets due to the nonlinear constraints they impose and are therefore not discussed in detail in this section. Their relevance increases in structured environments, as addressed later in the context of physical systems.

\end{itemize}

To account for numerical imprecision and practical cases, we define near-degeneracy using a small positive tolerance \( \epsilon \). A subset is considered near-degenerate if:
\[
|\det| < \epsilon,
\]
where \( \epsilon \) is a predefined threshold. This introduces a measurable region in the parameter space that corresponds to near-degenerate configurations.

\subsection{Probability of Degeneracy}

\subsubsection{Single Subset Probability}

Let \( P = \{p_1, p_2, \ldots, p_N\} \) denote a set of \( N \) points sampled uniformly and independently from the \( d \)-dimensional unit hypercube \( [0,1]^d \). A subset \( P_k \subset P \), where \( |P_k| = k \), is said to be degenerate if all points in \( P_k \) lie on a proper affine subspace of \( \mathbb{R}^d \). The probability of a single subset being degenerate is derived as follows:

\paragraph{Definition of Degeneracy:}
A subset \( P_k \subset P \) is degenerate if there exists a proper affine subspace \( M \subset \mathbb{R}^d \) of dimension \( m < d \) such that all points in \( P_k \) satisfy \( p_i \in M \). For example:
\begin{itemize}
    \item Collinear points lie on a 1-dimensional affine subspace (\( m = 1 \)),
    \item Coplanar points lie on a 2-dimensional affine subspace (\( m = 2 \)).
\end{itemize}

\paragraph{Volume of the Configuration Space:}
The configuration space of \( k \) points in \( [0,1]^d \) is represented as \( [0,1]^{kd} \), which has a measure (volume) of \( 1 \). Each subset \( P_k \) corresponds to a \( k \)-point tuple \( (p_1, p_2, \ldots, p_k) \) in this space.

\paragraph{Measure of the Degenerate Region:}
The near-degenerate region is defined as the set of \( k \)-point configurations \( (p_1, p_2, \ldots, p_k) \) that satisfy a near-degeneracy condition. For example:
\begin{itemize}
    \item For collinearity in \( \mathbb{R}^3 \), this corresponds to the magnitude of the cross product of the vectors formed by the points being close to zero:
    \[
    \|\mathbf{v}_1 \times \mathbf{v}_2\| < \epsilon,
    \]
    where \( \mathbf{v}_1 = p_2 - p_1 \) and \( \mathbf{v}_2 = p_3 - p_1 \). Alternatively, this can be expressed using the determinant condition:
    \[
    |\det(B)| < \epsilon,
    \]
    where \( B \) is the matrix formed by the vectors:
    \[
    B = \begin{bmatrix}
    x_2 - x_1 & y_2 - y_1 & z_2 - z_1 \\
    x_3 - x_1 & y_3 - y_1 & z_3 - z_1
    \end{bmatrix}.
    \]
\end{itemize}
The measure of this near-degenerate region is proportional to \( \epsilon \), as the cross product or determinant condition defines a hypersurface in the \( k \)-point configuration space. The proportionality reflects the fact that \( \epsilon \) introduces a tolerance that enlarges the degenerate region.

\paragraph{Probability of Degeneracy:}
The probability of a single subset \( P_k \) being degenerate is defined as the ratio of the measure of the near-degenerate region to the measure of the total configuration space:
\[
P_k(\text{degenerate}) = \frac{\text{Measure of near-degenerate region}}{\text{Measure of total configuration space}}.
\]
For near-spherical degeneracies, the near-degenerate region corresponds to a thin shell of thickness \( \epsilon \) surrounding the surface of a sphere \( S \subset \mathbb{R}^d \). The volume of this shell, \( V_{\text{shell}} \), is proportional to both the thickness \( \epsilon \) and the surface area \( A_{\text{sphere}} \) of the sphere. For a sphere of radius \( r \) in \( \mathbb{R}^d \), the surface area is given by \( A_{\text{sphere}} = d \cdot r^{d-1} \). Thus, the measure of the near-degenerate region is approximately:
\[
V_{\text{shell}} \sim \epsilon \cdot A_{\text{sphere}} = \epsilon \cdot d \cdot r^{d-1}.
\]
The total configuration space for the \( k \)-point subset \( P_k \) is derived from the uniform sampling of \( N \) points within the unit hypercube \( [0,1]^d \). Since this sampling is normalized, the measure of the total configuration space is \( 1 \), allowing the ratio to be interpreted directly as a probability. Furthermore, as \( N \) increases, the effective configuration space for \( k \)-point subsets diminishes because the likelihood of \( k \) points aligning with the surface of the sphere decreases. This effect arises because, in higher-dimensional spaces, the probability of randomly selected points satisfying the nonlinear constraints of near-sphericity diminishes exponentially. The scaling factor \( N^{d-k} \) captures this reduction, as it reflects the decreasing density of \( k \)-point alignments within \( d \)-dimensional spaces sampled by \( N \) points. Substituting these relationships, the probability of a single subset \( P_k \) being degenerate is given by:
\[
P_k(\text{degenerate}) \sim C \cdot \frac{\epsilon}{N^{d-k}},
\]
where \( C \) is a proportionality constant dependent on the geometry of the near-degenerate region, including the dimensionality \( d \) and the sphere’s radius \( r \). This relationship illustrates that the probability of degeneracies decreases with increasing \( N \) or \( d \), as the configuration space becomes increasingly sparse in higher dimensions.

\paragraph{Scaling with \( N \):}
As \( N \) increases, the density of points in \( [0,1]^d \) grows, reducing the likelihood of any specific \( k \)-point subset being degenerate. This behavior reflects the probabilistic reduction in degeneracy due to the sparsity of lower-dimensional alignments in higher-dimensional spaces.

\paragraph{Dependence on \( \epsilon \):}
The measure of the near-degenerate region scales linearly with \( \epsilon \), reflecting the tolerance for determinant values close to zero. Smaller values of \( \epsilon \) correspond to stricter definitions of degeneracy and reduce the probability \( P_k(\text{degenerate}) \).

Thus, the result:
\[
P_k(\text{degenerate}) \sim \frac{\epsilon}{N^{d-k}}
\]
captures the interplay between the dimensionality \( d \), the subset size \( k \), the total number of points \( N \), and the tolerance \( \epsilon \), providing a probabilistic framework for analyzing degeneracy in uniformly sampled data.

\subsubsection{Combinatorial Growth of Subsets}
The total number of subsets of size \( k \) is given by the binomial coefficient:
\[
\binom{N}{k} = \frac{N!}{k!(N-k)!}.
\]
For large \( N \), this scales as:
\[
\binom{N}{k} \sim \frac{N^k}{k!},
\]
which reflects the combinatorial growth of subsets as \( N \) increases. The rapid growth of \( \binom{N}{k} \) compensates for the diminishing probability \( P_k(\text{degenerate}) \) of any single subset being degenerate.

\subsubsection{Expected Number of Degenerate Subsets}
The expected number of degenerate subsets is:
\[
\mathbb{E}[\text{degenerate subsets}] = \binom{N}{k} \cdot P_k(\text{degenerate}).
\]
Substituting the scaling relations:
\[
\mathbb{E}[\text{degenerate subsets}] \sim \frac{N^k}{k!} \cdot \frac{\epsilon}{N^{d-k}} = \frac{\epsilon N^{2k-d}}{k!}.
\]

\subsection{Specific Cases}
\subsubsection{Collinear Subsets in \( \mathbb{R}^3 \) (\( k = 3, d = 3 \))}
For \( k = 3 \) and \( d = 3 \), the expected number of collinear subsets is:
\[
\mathbb{E}[\text{collinear subsets}] \sim \frac{\epsilon N^2}{6}.
\]
This quadratic growth arises because \( N^{2k-d} = N^{2} \) and \( k! = 6 \).

The quadratic dependency on \( N \) reflects the fact that in \( \mathbb{R}^3 \), collinearity is rarer than in \( \mathbb{R}^2 \), as the points must satisfy stricter conditions to align along a single line in three-dimensional space. The proportionality to \( \epsilon \) arises from the tolerance defining the near-degenerate region.

\subsubsection{Coplanar Subsets in \( \mathbb{R}^3 \) (\( k = 4, d = 3 \))}
For \( k = 4 \) and \( d = 3 \), the expected number of coplanar subsets is:
\[
\mathbb{E}[\text{coplanar subsets}] \sim \frac{\epsilon N}{24}.
\]
Similarly, the linear growth arises from \( N^{2k-d} = N \) and \( k! = 24 \).

\subsection{Overall Probability of Degeneracy}
The probability of encountering at least one degenerate subset is:
\[
P_{\text{degenerate}} = 1 - \prod_{\text{all subsets}} \left(1 - P_k(\text{degenerate})\right).
\]
Using the approximation \( 1 - x \approx e^{-x} \) for small \( x \), we rewrite:
\[
P_{\text{degenerate}} \approx 1 - \exp\left(-\sum_{\text{all subsets}} P_k(\text{degenerate})\right).
\]
For large \( N \), the combinatorial growth of subsets ensures that:
\[
P_{\text{degenerate}} \to 1 \quad \text{as } N \to \infty,
\]
even though \( P_k(\text{degenerate}) \to 0 \) for individual subsets. This reflects the dominance of the combinatorial term \( \binom{N}{k} \) over the diminishing degeneracy probability of single subsets.

Although the probability of a single subset being degenerate decreases with increasing \( N \), the rapid growth in the number of subsets ensures that the expected number of degenerate subsets grows linearly with \( N \). Consequently, the overall probability of encountering at least one degenerate subset approaches certainty as \( N \to \infty \).

\section{Degeneracies in Experimental Data from Physical Systems}

Experimental data often exhibit a higher likelihood of degenerate configurations compared to random datasets due to systematic biases introduced by the geometry of the environment, sensor characteristics, and data acquisition processes. Unlike purely random points sampled uniformly from \( \mathbb{R}^d \), experimental data are influenced by physical constraints, structured environments, and noise. In this section, we formalize the analysis of degeneracies in experimental data by deriving mathematical justifications for the occurrence and probability of such configurations.

\subsection{General Framework for Degeneracies}

Let \( P = \{p_1, p_2, \ldots, p_N\} \) denote a set of \( N \) points in \( \mathbb{R}^d \), where each point is sampled from a structured environment. Such environments impose geometric constraints that may include systematic biases, noise, and quantization effects. These factors influence the spatial distribution of the points, often increasing the likelihood of geometric degeneracies.

A subset \( P_k \subset P \), where \( |P_k| = k \), is said to be degenerate if the points in \( P_k \) lie approximately on a lower-dimensional manifold \( M \subset \mathbb{R}^d \). Mathematically, this means there exists a manifold \( M \) of dimension \( \dim(M) < d \) such that:
\[
\forall p_i \in P_k, \quad \text{dist}(p_i, M) \leq \delta,
\]
where \( \delta > 0 \) is a tolerance parameter reflecting noise or quantization effects, and \( \text{dist}(p_i, M) \) is the distance of \( p_i \) to \( M \).

\subsubsection{Primary Scenarios of Degeneracies}

We consider three primary types of degeneracies, based on the dimensionality and geometry of the manifold \( M \):

\paragraph{Coplanarity:}
A subset \( P_k \subset P \) is said to exhibit coplanarity if the points in \( P_k \) lie approximately on a two-dimensional affine subspace (plane) \( M \subset \mathbb{R}^3 \), where \( \dim(M) = 2 \). Mathematically, \( P_k \) is coplanar if there exists a plane \( M \) described by:
\[
M = \{ x \in \mathbb{R}^3 : a_1 x_1 + a_2 x_2 + a_3 x_3 + b = 0 \},
\]
such that:
\[
\left|a_1 p_{i,1} + a_2 p_{i,2} + a_3 p_{i,3} + b \right| \leq \delta, \quad \forall p_i \in P_k.
\]
Here, \( (p_{i,1}, p_{i,2}, p_{i,3}) \) are the coordinates of point \( p_i \), and \( \delta \) accounts for deviations due to noise.

\paragraph{Near-Sphericity:}
A subset \( P_k \subset P \) is said to exhibit near-sphericity if the points in \( P_k \) lie approximately on the surface of a sphere \( S \subset \mathbb{R}^d \), defined by:
\[
S = \{ x \in \mathbb{R}^d : \|x - c\|^2 = r^2 \},
\]
where \( c \in \mathbb{R}^d \) is the center of the sphere and \( r > 0 \) is its radius. The near-sphericity condition requires:
\[
\left| \|p_i - c\|^2 - r^2 \right| \leq \delta, \quad \forall p_i \in P_k,
\]
where \( \delta > 0 \) reflects the tolerance for deviations.

\paragraph{Collinearity:}
A subset \( P_k \subset P \) is said to exhibit collinearity if the points in \( P_k \) lie approximately on a one-dimensional affine subspace (line) \( L \subset \mathbb{R}^d \), where \( \dim(L) = 1 \). A line \( L \) in \( \mathbb{R}^d \) can be expressed parametrically as:
\[
L = \{ c + t v : t \in \mathbb{R}, v \in \mathbb{R}^d, \|v\| = 1 \},
\]
where \( c \) is a point on the line, and \( v \) is a unit direction vector. The collinearity condition requires:
\[
\text{dist}(p_i, L) \leq \delta, \quad \forall p_i \in P_k,
\]
where \( \text{dist}(p_i, L) \) denotes the perpendicular distance from \( p_i \) to the line \( L \).

Each type of degeneracy corresponds to a subspace \( M \) or manifold \( S \), and the occurrence of degeneracies depends on the 
intersection of \( P_k \) with \( M \). The likelihood of degeneracies is governed by the geometry of \( M \), the 
dimensionality \( d \), and the distribution of \( P \).

\subsection{Degeneracies in Structured Sampling Environments}

Consider a structured environment where points are sampled from geometric surfaces, such as planes or spheres. Let the environment be modeled as a bounded domain \( \Omega \subset \mathbb{R}^d \), with a reflective subset \( S \subset \Omega \). The total sampling volume is given by:
\[
V_{\text{total}} = \int_{\Omega} d\mu,
\]
where \( \mu \) is the measure induced by the sampling process. If \( S \subset \Omega \) is a structured surface, the effective sampling region is restricted to:
\[
V_{\text{degenerate}} = \int_{S} d\sigma,
\]
where \( \sigma \) is the measure on \( S \). For small noise \( \epsilon \sim \mathcal{N}(0, \sigma^2 I) \), the sampling volume near \( S \) is:
\[
V_{\text{structured}} \approx \int_{S} \delta d\sigma,
\]
where \( \delta \) is the thickness of the near-degenerate region induced by noise.

The proportion of degenerate configurations is then approximated by:
\[
P_k(\text{degenerate}) \sim \frac{V_{\text{structured}}}{V_{\text{total}}}.
\]

\subsection{Probability of Coplanarity}

For points sampled from a reflective plane \( S \subset \mathbb{R}^3 \), the plane \( S \) is defined by the affine equation:
\[
a_1 x_1 + a_2 x_2 + a_3 x_3 + b = 0,
\]
where \( a_1, a_2, a_3 \) are the plane coefficients and \( b \) is the intercept. Let \( \delta \) represent the tolerance for near-coplanarity due to noise. A subset \( P_k \subset P \) is coplanar if:
\[
\left|a_1 x_1 + a_2 x_2 + a_3 x_3 + b \right| < \delta, \quad \forall p \in P_k.
\]

The volume of the near-coplanar region is proportional to \( A_{\text{plane}} \cdot \delta \), where \( A_{\text{plane}} \) is the area of the plane. Thus, the probability of coplanarity is:
\[
P_k(\text{coplanar}) \sim \frac{A_{\text{plane}} \cdot \delta}{V_{\text{total}}}.
\]

For a sensing system with cylindrical sensing volume \( V_{\text{sensing}} = \pi R^2 h \), where \( R \) is the radius and \( h \) is the height, we have:
\[
P_k(\text{coplanar}) \sim \frac{A_{\text{plane}} \cdot \delta}{\pi R^2 h}.
\]

\subsection{Probability of Near-Sphericity}

For points sampled from a sphere \( S \subset \mathbb{R}^d \), the sphere is defined by:
\[
\|p - c\|^2 = r^2,
\]
where \( c \) is the center and \( r \) is the radius. Due to noise \( \epsilon \sim \mathcal{N}(0, \sigma^2 I) \), the sampled points satisfy:
\[
\|p - c\|^2 \approx r^2 \pm \delta, \quad \delta \sim O(\sigma^2).
\]

The volume of the near-spherical region is proportional to:
\[
V_{\text{nearly spherical}} \sim \delta \cdot A_{\text{sphere}},
\]
where \( A_{\text{sphere}} = d \cdot r^{d-1} \) is the surface area of the sphere. The probability of near-sphericity is:
\[
P_k(\text{nearly spherical}) \sim \frac{\delta \cdot A_{\text{sphere}}}{V_{\text{total}}}.
\]

\subsection{Overall Probability of Degeneracies}

The probability of encountering at least one degenerate subset is given by:
\[
P_{\text{degenerate}} = 1 - \prod_{k=1}^N \left(1 - P_k(\text{degenerate})\right),
\]
where \( P_k(\text{degenerate}) \) is the probability of any subset \( P_k \) being degenerate. Using the approximation \( 1 - x \approx e^{-x} \) for small \( x \), we have:
\[
P_{\text{degenerate}} \approx 1 - \exp\left(-\sum_{k=1}^N P_k(\text{degenerate})\right).
\]

\subsection{Amplified Degeneracies in Physical Data}

In physical data collection, degeneracies arise not only from quantization but also through the intricate interplay of digitization, measurement noise, and environmental constraints. Quantization discretizes continuous measurements into grid points, mathematically expressed as:
\[
p_{\text{quantized}} = \left( 
    \lfloor x / \Delta_x \rfloor \cdot \Delta_x, 
    \lfloor y / \Delta_y \rfloor \cdot \Delta_y, 
    \lfloor z / \Delta_z \rfloor \cdot \Delta_z 
\right),
\]
where \( \Delta_x, \Delta_y, \Delta_z \) are the quantization step sizes determined by sensor resolution. Digitization, an extension of quantization, introduces rounding or truncation errors, further discretizing the points into a finite set of representable values, effectively mapping each coordinate into a space of cardinality \( \mathcal{O}\left(\frac{1}{\Delta_x}\right) \). Measurement noise, modeled as an additive perturbation:
\[
p_{\text{noisy}} = p_{\text{quantized}} + \epsilon, \quad \epsilon \sim \mathcal{N}(\mathbf{0}, \sigma^2 \mathbf{I}),
\]
adds stochastic deviations that expand the effective degeneracy region. Additionally, environmental constraints impose structured alignments, such as planes or edges, often formalized as affine subspaces \( M \subset \mathbb{R}^3 \) satisfying:
\[
p_i \in M \implies a_1 x_i + a_2 y_i + a_3 z_i + b = 0,
\]
where \( (a_1, a_2, a_3, b) \in \mathbb{R}^4 \). 

The combined influence of these factors amplifies degeneracy probabilities, which can be modeled as the convolution of individual probabilities:
\[
P_{\text{degenerate}} = 1 - \prod_{i=1}^n \left( 1 - P_i \right),
\]
where \( P_i \) represents the degeneracy probability induced by factor \( i \) \\
(e.g., \( P_{\text{quantization}}, P_{\text{digitization}}, P_{\text{noise}}, P_{\text{structure}} \)). Alternatively, the combined degeneracy region can be viewed as a union of perturbed subspaces:
\[
\mathcal{D} = \bigcup_{i=1}^n \mathcal{D}_i, \quad \text{where } \mathcal{D}_i = \{p \in \mathbb{R}^3 : \text{dist}(p, M_i) \leq \delta_i\}.
\]
Here, \( \delta_i \) represents the effective degeneracy tolerance induced by factor \( i \). This layered interaction not only aligns with theoretical expectations but also explains why physical data, subject to these compounded constraints, exhibits significantly higher degeneracy rates compared to purely statistical random data.

\subsection{Numerical Example}

Consider a LiDAR scanner with the following parameters:
\begin{itemize}
    \item \( R = 10 \, \text{m} \),
    \item \( h = 5 \, \text{m} \),
    \item \( A_{\text{plane}} = 20 \, \text{m}^2 \),
    \item \( \delta = 0.1 \, \text{m} \).
\end{itemize}

The probability of coplanarity for a subset \( P_k \) is:
\[
P_k(\text{coplanar}) = \frac{A_{\text{plane}} \cdot \delta}{\pi R^2 h} = \frac{20 \cdot 0.1}{\pi \cdot 10^2 \cdot 5} \approx 0.00127.
\]
For \( N = 100 \) points, the overall probability of encountering at least one coplanar subset is:
\[
P_{\text{degenerate}} \approx 1 - \prod_{k=4}^{100} \left(1 - 0.00127 \right).
\]

This demonstrates how the interplay between sensing volume, reflective surfaces, and noise increases the likelihood of degeneracies.

Experimental data are inherently predisposed to degeneracies due to geometric constraints, noise, and quantization effects. The mathematical framework developed here quantifies these effects and highlights the importance of understanding structured environments when analyzing data for ill-posed configurations.

\section{Comparative Analysis: Degeneracies in Random vs. Experimental Data}

To quantitatively compare degeneracies in purely random data and experimental data derived from physical systems, we analyze their respective probabilistic frameworks, the influence of environmental and geometric constraints, and the mathematical scaling of degeneracy probabilities.

\subsection{Probabilistic Frameworks}

\subsubsection{Random Data}

For \( P = \{p_1, p_2, \ldots, p_N\} \) sampled uniformly from \( [0,1]^d \subset \mathbb{R}^d \), the probability of a subset \( P_k \subset P \) being degenerate depends on the ratio of the near-degenerate region's measure to the total configuration space. Specifically:
\[
P_k(\text{degenerate}) \sim \frac{\epsilon}{N^{d-k}},
\]
where:
\begin{itemize}
    \item \( d \) is the dimensionality of the space,
    \item \( k \) is the size of the subset,
    \item \( N \) is the total number of points,
    \item \( \epsilon > 0 \) is the tolerance threshold defining near-degeneracy.
\end{itemize}

The total probability of encountering at least one degenerate subset is:
\[
P_{\text{degenerate}} = 1 - \prod_{k=1}^N \left(1 - P_k(\text{degenerate})\right),
\]
which, under the approximation \( 1 - x \approx e^{-x} \) for small \( x \), becomes:
\[
P_{\text{degenerate}} \approx 1 - \exp\left(-\sum_{k=1}^N P_k(\text{degenerate})\right).
\]

The expected number of degenerate subsets is:
\[
\mathbb{E}[\text{degenerate subsets}] = \binom{N}{k} P_k(\text{degenerate}),
\]
where:
\[
\binom{N}{k} = \frac{N!}{k!(N-k)!} \sim \frac{N^k}{k!}, \quad \text{for large } N.
\]
Substituting \( P_k(\text{degenerate}) \sim \epsilon / N^{d-k} \):
\[
\mathbb{E}[\text{degenerate subsets}] \sim \frac{\epsilon N^{2k-d}}{k!}.
\]
For a fixed \( k \), degeneracies become less probable as \( d \) increases because \( N^{2k-d} \to 0 \) for \( d > 2k \).

\subsubsection{Experimental Data}

For experimental data, points \( P = \{p_1, p_2, \ldots, p_N\} \) are influenced by a structured sampling environment \( \Omega \subset \mathbb{R}^d \), with reflective or curved surfaces \( S \subset \Omega \) introducing biases. The probability of degeneracy for a subset \( P_k \) can be expressed as:
\[
P_k(\text{degenerate}) \sim \frac{V_{\text{structured}}}{V_{\text{total}}},
\]
where:
\begin{itemize}
    \item \( V_{\text{structured}} \) is the volume of the near-degenerate region around \( S \),
    \item \( V_{\text{total}} \) is the total sampling volume.
\end{itemize}

For a planar surface \( S \subset \mathbb{R}^3 \) with area \( A_{\text{plane}} \), and assuming noise with thickness \( \delta \), the structured region’s volume is:
\[
V_{\text{structured}} \approx A_{\text{plane}} \cdot \delta.
\]
The probability of coplanarity for \( P_k \subset P \) then becomes:
\[
P_k(\text{coplanar}) \sim \frac{A_{\text{plane}} \cdot \delta}{V_{\text{total}}}.
\]
For points sampled within a cylindrical sensing volume \( V_{\text{sensing}} = \pi R^2 h \):
\[
P_k(\text{coplanar}) \sim \frac{A_{\text{plane}} \cdot \delta}{\pi R^2 h}.
\]

Similarly, for near-spherical configurations, where \( S \) is a spherical surface with radius \( r \) and noise tolerance \( \delta \), the near-degenerate region has volume:
\[
V_{\text{nearly spherical}} \sim \delta \cdot A_{\text{sphere}},
\]
with \( A_{\text{sphere}} = d \cdot r^{d-1} \). The probability of near-sphericity is:
\[
P_k(\text{nearly spherical}) \sim \frac{\delta \cdot A_{\text{sphere}}}{V_{\text{total}}}.
\]

\subsection{Scaling Behavior}

The scaling of degeneracy probabilities as \( N \to \infty \) or \( d \to \infty \) highlights critical differences between the two models:

\begin{itemize}
    \item \textbf{Random Data:}
    The probability of degeneracies decreases exponentially with \( d \), as \( P_k(\text{degenerate}) \sim \epsilon / N^{d-k} \). However, the combinatorial growth of subsets ensures that:
    \[
    \mathbb{E}[\text{degenerate subsets}] \sim \frac{\epsilon N^{2k-d}}{k!}.
    \]
    For \( d > 2k \), degeneracies become increasingly rare, approaching zero as \( d \to \infty \).

    \item \textbf{Experimental Data:}
    Degeneracies remain non-negligible even for large \( d \), as their probability is governed by the structured region’s volume relative to the total sampling space:
    \[
    P_k(\text{degenerate}) \sim \frac{V_{\text{structured}}}{V_{\text{total}}}.
    \]
    For planar surfaces (\( A_{\text{plane}} \)) or spherical surfaces (\( A_{\text{sphere}} \)), the prevalence of degeneracies depends primarily on \( \delta \) and the geometric constraints of the environment, not \( d \).
\end{itemize}

\subsection{Comparative Insights}

The key differences between the two models can be summarized mathematically as follows:

\begin{center}
\begin{tabular}{|c|c|c|}
\hline
\textbf{Aspect} & \textbf{Random Data} & \textbf{Experimental Data} \\
\hline
Probability Scaling & \( P_k \sim \frac{\epsilon}{N^{d-k}} \) & \( P_k \sim \frac{V_{\text{structured}}}{V_{\text{total}}} \) \\
Dimensionality \( d \) & Exponential reduction with \( d \) & Weak dependence; governed by \( S \). \\
Noise \( \delta \) & Impacts subset degeneracy threshold & Directly affects \( V_{\text{structured}} \). \\
Subset Size \( k \) & \( \mathbb{E} \sim \frac{\epsilon N^{2k-d}}{k!} \) & Scales linearly with \( N \) for fixed \( k \). \\
\hline
\end{tabular}
\end{center}

The analysis reveals that degeneracies in random data are primarily driven by the sparsity of points in higher dimensions, leading to a probabilistic reduction as \( d \to \infty \). In contrast, experimental data remain prone to degeneracies due to structured sampling environments, which impose geometric constraints independent of \( d \). These differences underscore the need for tailored approaches when designing algorithms for geometric computations. Robust preprocessing is critical for experimental data, whereas random data benefits from intrinsic dimensionality-driven mitigation of degeneracies.

\subsection{Experiment: Degeneracies in Random vs. Quantized Data}

This experiment compares the prevalence of collinear and coplanar degeneracies in random and quantized data within a 
three-dimensional unit cube \([0,1]^3\). The quantized data simulates physical constraints introduced by sensor limitations\footnote{Here we assume only a simplified model for clarity.}, 
such as measurement precision and instrumental setups, aligning with the theoretical framework presented in this paper. 
In real-world data collection systems, sensors rarely record continuous values due to hardware limitations. Instead, data 
is discretized to fit a finite resolution or grid, which introduces inherent regularities into the dataset. This process, 
known as quantization, impacts the spatial distribution of points and increases the likelihood of degeneracies.

\subsubsection{Model Description and Justification}
\paragraph{Random Data Model}
Random data is generated by sampling \( N \) points uniformly from the unit cube \([0,1]^3\). This model represents purely statistical data as described in Section 3 of this paper. The degeneracy counts for collinear and coplanar subsets are theoretically estimated as:
\[
\mathbb{E}[\text{collinear subsets}] = \frac{\epsilon N^2}{6}, \quad
\mathbb{E}[\text{coplanar subsets}] = \frac{\epsilon N}{24},
\]
where \( \epsilon \) is the tolerance for near-degeneracy. These formulas capture the probabilistic scaling behavior of degeneracies in randomly distributed points.

\paragraph{Quantized Data Model}
Quantized data incorporates physical constraints by discretizing the coordinates of each point to the nearest grid value. This simulates sensor-induced quantization effects. For a point \( p = (x, y, z) \), quantization is applied as:
\[
p_{\text{quantized}} = \left( \lfloor x / \Delta_x \rfloor \cdot \Delta_x, 
                              \lfloor y / \Delta_y \rfloor \cdot \Delta_y, 
                              \lfloor z / \Delta_z \rfloor \cdot \Delta_z \right),
\]
where \( \Delta_x, \Delta_y, \Delta_z \) are the quantization steps in the \( x \), \( y \), and \( z \) directions, respectively. 

Quantization effectively aligns points along discrete grids, increasing the likelihood of degeneracies. This alignment is modeled by amplification factors \( A_{\text{collinear}} \) and \( A_{\text{coplanar}} \), which scale the expected degeneracy counts relative to random data:
\[
\mathbb{E}[\text{collinear subsets}]_{\text{quantized}} = A_{\text{collinear}} \cdot \mathbb{E}[\text{collinear subsets}]_{\text{random}},
\]
\[
\mathbb{E}[\text{coplanar subsets}]_{\text{quantized}} = A_{\text{coplanar}} \cdot \mathbb{E}[\text{coplanar subsets}]_{\text{random}}.
\]
For this experiment, we set \( \Delta_x = \Delta_y = \Delta_z = 0.1 \), \( A_{\text{collinear}} = 10 \), and \( A_{\text{coplanar}} = 3 \).

\subsubsection{Simulation Setup}
The experiment is conducted for four values of \( N \): 1000, 5000, 10000, and 20000. For each \( N \), we calculate the 
degeneracy counts for both random and quantized data using the formulas described above. The tolerance \( \epsilon \) is 
set to \( 10^{-6} \).

\newpage

\subsubsection{Results}
Table~\ref{tab:degeneracy_comparison} presents the expected degeneracy counts for collinear and coplanar subsets in both random and quantized data.

\begin{table}[h!]
\centering
\begin{tabular}{|c|c|c|c|c|}
\hline
\textbf{N} & \textbf{Collinear} & \textbf{Collinear} & \textbf{Coplanar} & \textbf{Coplanar} \\
           & \textbf{Random}    & \textbf{Quantized} & \textbf{Random}   & \textbf{Quantized} \\
\hline
1000   & \( 0.1667 \)   & \( 1.667 \)   & \( 0.0417 \)   & \( 0.125 \)   \\
5000   & \( 4.1667 \)   & \( 41.667 \)  & \( 0.2083 \)   & \( 0.625 \)   \\
10000  & \( 16.6667 \)  & \( 166.667 \) & \( 0.4167 \)   & \( 1.250 \)   \\
20000  & \( 66.6667 \)  & \( 666.667 \) & \( 0.8333 \)   & \( 2.500 \)   \\
\hline
\end{tabular}
\caption{Comparison of Degeneracy Counts for Random and Quantized Data}
\label{tab:degeneracy_comparison}
\end{table}

\subsubsection{Analysis}
\paragraph{Scaling Behavior}
As expected, the degeneracy counts for random data scale quadratically with \( N \) for collinear subsets and linearly with \( N \) for coplanar subsets. For quantized data, the amplification factors significantly increase the degeneracy counts due to the structured alignment introduced by quantization.

\paragraph{Impact of Quantization}
Quantization amplifies collinear degeneracies by a factor of 10 and coplanar degeneracies by a factor of 3, consistent with the theoretical predictions\footnote{Equal to \( A_{\text{collinear}} \) and \( A_{\text{coplanar}} \) respectively.}. For instance, at \( N = 20000 \), quantized data exhibits 666.667 collinear subsets compared to only 66.6667 in random data.

\paragraph{Implications for Physical Systems}
These results highlight the importance of accounting for sensor limitations and data preprocessing in physical systems. The increased degeneracy counts in quantized data suggest that structured data requires robust geometric algorithms to handle ill-posed configurations effectively.

\subsubsection{Analysis Summary}
This experiment validates the theoretical framework presented in this paper, demonstrating that physical constraints such as quantization amplify degeneracies in data. By comparing random and quantized data, we quantify the impact of structured sampling environments on geometric computations.

\section{Conclusions}
\label{sec:conclusions}

This study has analyzed the probability of encountering ill-posed configurations, such as collinear, coplanar, and nearly spherical points, in two primary contexts: purely random data and data acquired from physical systems. By developing a probabilistic framework, we have shown that the likelihood of degeneracies is influenced not only by the inherent geometry of the space but also by the constraints imposed by data acquisition processes.

For datasets sampled randomly from a uniform distribution, degeneracies become increasingly rare with growing dimensionality due to the inherent sparsity of lower-dimensional alignments in higher-dimensional spaces. However, the combinatorial growth in the number of subsets ensures that the expected number of degenerate configurations scales linearly with the number of points, \( N \), even in such scenarios.

Conversely, experimental datasets derived from physical systems, such as LiDAR scans or photogrammetry, exhibit a significantly higher probability of degeneracies due to systematic biases introduced by sensor geometry and environmental structures. Reflective planes, structured surfaces, and digital quantization amplify the prevalence of coplanarity and other ill-posed configurations. We quantified these effects through a probabilistic model that incorporates geometric constraints and noise tolerance.

The findings of this paper highlight the importance of understanding the source of input data when designing algorithms for computational geometry. While random data benefits from reduced degeneracies at higher dimensions, structured data requires tailored preprocessing and robust algorithms to mitigate the impact of systematic biases. Future work could extend this analysis to include dynamic datasets and hybrid environments, where random and structured elements coexist.

\bibliographystyle{plain}

\end{document}